\documentclass[sn-aps]{sn-jnl}

\usepackage[utf8]{inputenc}
\usepackage[T1]{fontenc}
\usepackage{amsmath,stackengine}
\usepackage{amsfonts}
\usepackage{amssymb}
\usepackage{listings}
\usepackage{hyperref}
\usepackage{graphicx}
\usepackage{xcolor}
\usepackage{enumerate}
\usepackage{mathtools}
\usepackage{subcaption}
\usepackage{manyfoot}%
\usepackage[linesnumbered,ruled,vlined]{algorithm2e}
\usepackage{grffile}
\usepackage{url}
\usepackage{stfloats}
\usepackage{textcomp}

\begin{document}
\title[LQT state feedback]{Data-driven Linear Quadratic Tracking based Temperature Control of a Big Area Additive Manufacturing System}

\author*[1,2,3]{\fnm{Eleni} \sur{Zavrakli}}\email{Eleni.Zavrakli@mu.ie}

\author[1,2,3]{\fnm{Andrew} \sur{Parnell}}

\author[3,5]{\fnm{Andrew} \sur{Dickson}}

\author[4]{\fnm{Subhrakanti} \sur{Dey}}

\affil*[1]{\orgdiv{Hamilton Institute}, \orgname{Maynooth University},  \state{Co. Kildare}, \country{Ireland}}

\affil[2]{\orgdiv{Department of Mathematics and Statistics}, \orgname{Maynooth University},  \state{Co. Kildare}, \country{Ireland}}

\affil[3]{\orgdiv{I-Form}, \orgname{Advanced Manufacturing Research Centre}, \country{Ireland}}

\affil[4]{\orgdiv{Division of Signals and Systems, Department of Electrical Engineering}, \orgname{Uppsala University}, \country{Sweden}}
\affil[5]{\orgdiv{School of Mechanical \& Materials Engineering}, \orgname{University College Dublin}, \country{Ireland}}

\abstract{Designing efficient closed-loop control algorithms is a key issue in Additive Manufacturing (AM), as various aspects of the AM process require continuous monitoring and regulation, with temperature being a particularly significant factor. Here we study closed-loop control of a state space temperature model with a focus on both model-based and data-driven methods. We demonstrate these approaches using a simulator of the temperature evolution in the extruder of a Big Area Additive Manufacturing system (BAAM). We perform an in-depth comparison of the performance of these methods using the simulator. We find that we can learn an effective controller using solely simulated process data. Our approach achieves parity in performance compared to model-based controllers and so lessens the need for estimating a large number of parameters of the intricate and complicated process model. We believe this result is an important step towards autonomous intelligent manufacturing. }




\maketitle

\section{Introduction}

Additive Manufacturing (AM) is shaping up to be one of the most promising methods in manufacturing, with its application ranging from medicine, aerospace, architecture to most aspects of daily life \cite{DILBEROGLU2017545,BERMAN2012155,BUCHANAN2019332}. AM has become popular partially because of its ability to handle very complex geometries and due to bringing manufacturing into the hands of the consumer. However, while AM processes are becoming more widely adopted, there is still a substantial lack of efficient closed-loop control algorithms \cite{mercado2020additive}. More specifically, there are many aspects of the process that require constant monitoring and control, some specific to a type of AM, others more universal to manufacturing systems. These aspects include temperature, printing speed, layer thickness, material viscosity and flow, and extrusion pressure, among others. 

Control theory \cite{anderson2007optimal,ogata2010modern} offers an array of methods that address the problem of controlling a system's behaviour. Most classic approaches involve designing controllers using a model of the system dynamics. This presents a very important limitation, as an exact model is not always available. The lack of knowledge of the system dynamics gave rise to the area of data-driven control \cite{hou2013model,de2019formulas,rosolia2017learning,piga2017direct}. Reinforcement Learning (RL) \cite{si2004handbook,sutton2018reinforcement} is a family of data-driven Machine Learning algorithms that address the problem of finding optimal policies within complex and uncertain environments through maximising rewards (or minimising costs). In the Control Theory literature, optimal policies are created via Dynamic Programming \cite{bellman1966dynamic,bertsekas2012dynamic,bertsekas1996neuro} with many similarities to RL but with differing terminology. Whilst some RL approaches have been developed in a deterministic fashion assuming knowledge of the system dynamics, most RL applications are data-driven and model-free. 

Some previous approaches have applied control algorithms to address specific challenges in metal AM, including the design of Proportional Integral Derivative (PID) controllers for the cooling rate \cite{farshidianfar2016real-time} and laser power \cite{HU200351} in laser-based AM,
the meltpool in Selective Laser Melting \cite{kruth2007feedback}, and layer height for Laser Metal Deposition \cite{8301604}. In recent years, with RL gaining in popularity, researchers have applied it to solve various AM-related problems including improving process monitoring through imaging \cite{8362941}, acoustic emissions \cite{wasmer2019situ}, and solving the machine scheduling problem \cite{alicastro2021reinforcement}. Additionally, RL has been used in toolpath generation  \cite{patrick2018reinforcement} and optimising melt pool depth \cite{ogoke2021thermal}. An RL framework has been used in a model-based application for implementing corrective actions on a robot wire arm AM system \cite{dharmawan2020model}. 

We study the problem of designing closed-loop controllers for a state space temperature model. State space models are very useful in the design of optimal controllers \cite{anderson2007optimal,ogata2010modern} but have not been widely used in AM applications. Previous work in state space modelling in manufacturing includes the identification of a state space model and its use to design a controller for temperature regulation \cite{gootjes2017applying}, the development of a state space model for product quality \cite{stoyanov2017machine}, and the use of a state space model for observer design of temperature states within an AM-produced part \cite{8430941}. 

The control objective throughout our paper can be considered as the optimisation of a quadratic performance index, namely the Linear Quadratic Tracking (LQT) problem. Our paper is structured as follows. Section \ref{baam} provides a brief introduction into the family of manufacturing systems that are the subject of this work. Section \ref{section_model} introduces the problem statement and the particular system that is used as a case study. We introduce a model describing the temperatures within the extruder head of a Big Area Additive Manufacturing (BAAM) system in state space form. In Section \ref{section model based}, we perform an in-depth analysis of model-based control methods. We explore the finite horizon LQT controller, and for the infinite-horizon case, we create a model-based RL controller. The model-based methods designed in this section can be applied to any state space model. In Section \ref{section q-learning}, we develop a model-free, data-driven controller. We use a discrete time, continuous state and action space Q-learning algorithm. Finally, in Section \ref{section results}, we present the setup of our simulations, followed by a detailed comparison of the performance of all the methods. Apart from simulating the system's behaviour after each controller is implemented, the system model is used for the design of all model-based controllers, and for generating the data that will be used to train the data-driven controller. We provide a thorough discussion on the results in Section \ref{section_conclusion}.

\section{BAAM - Big Area Additive Manufacturing heating systems}\label{baam}

BAAM style 3D printing systems are becoming more widely utilised in the areas of large scale polymer and composite manufacturing. This technology is a derivation of the MEX (Material Extrusion) technology prevalent in desktop sized 3D printers (Examples include Prusa i3 and Ultimaker S5). Unlike the majority of desktop MEX systems which contain one heating source (Typically a ceramic heater cartridge) and a sensor (Thermistor, Thermocouple or similar), BAAM systems contain multiple heaters at different points along the length of their heated melt zone. The primary reason for this is the substantially larger volume of material that is processed through a BAAM extruder per unit of time, which can range from 2 - 20Kg per hour, compared with desktop systems which would average 0.05Kg per hour. As polymers are almost exclusively insulative materials (although sometimes containing conductive additives such as Carbon black or metallic powders), heat transfer from the melt zone wall into the polymer itself is poor, and thus it is necessary to increase the length of time the polymer resides inside the melt zone so as to evenly melt the polymer granules. In order to reach the desired residency times whilst maintaining a high throughput of material, melt zones are staggered across a long heated bore (which can commonly range from 0.2 to 0.5 meters in length). 

\begin{figure}[ht]
    \centering
    \includegraphics[width=\textwidth,height=7cm]{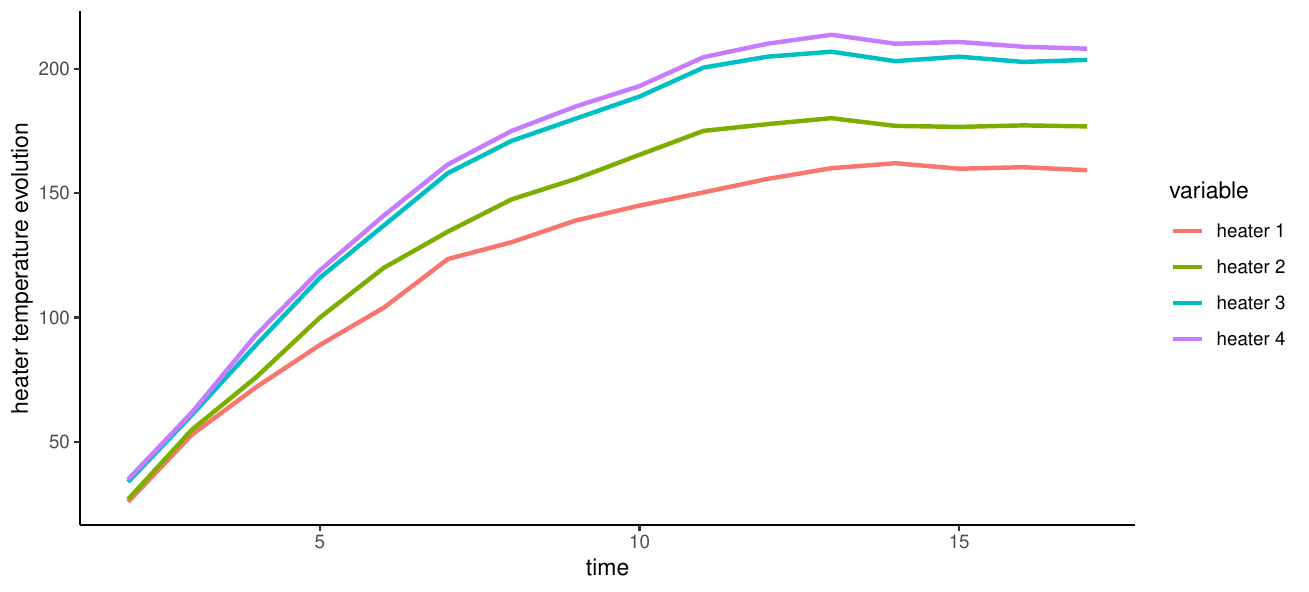}
    \caption{Example of heater staggering for a BAAM system using PLA, consisting of 4 heating zones to melt the material.}
    \label{temp example}
\end{figure}

The benefit of a multi-zone heating system is the accurate control of the melt pool, which is typically arranged in a steadily increasing temperature so as to slowly increase polymer temperature until the optimum state is reached for extrusion. An example of this heating system is visually represented in Figure \ref{temp example}. It is therefore vital that a robust method of controlling these temperatures is used, particularly as this technology is adopted for printing of PEEK, PEI and other high value polymers, where the cost of a print failure can be significant. Finally, in order to ensure that a temperature is maintained throughout the printing process, it is essential that an adaptive technique be adopted to ensure that temperature variance is minimal, even as unavoidable systems changes such as flow rate, ambient temperature, power supply etc occur during the print.

\section{Linear Quadratic Tracking Control of a State-Space Temperature System} \label{section_model}

\begin{figure}[ht]
    \centering
    \includegraphics[width=\textwidth]{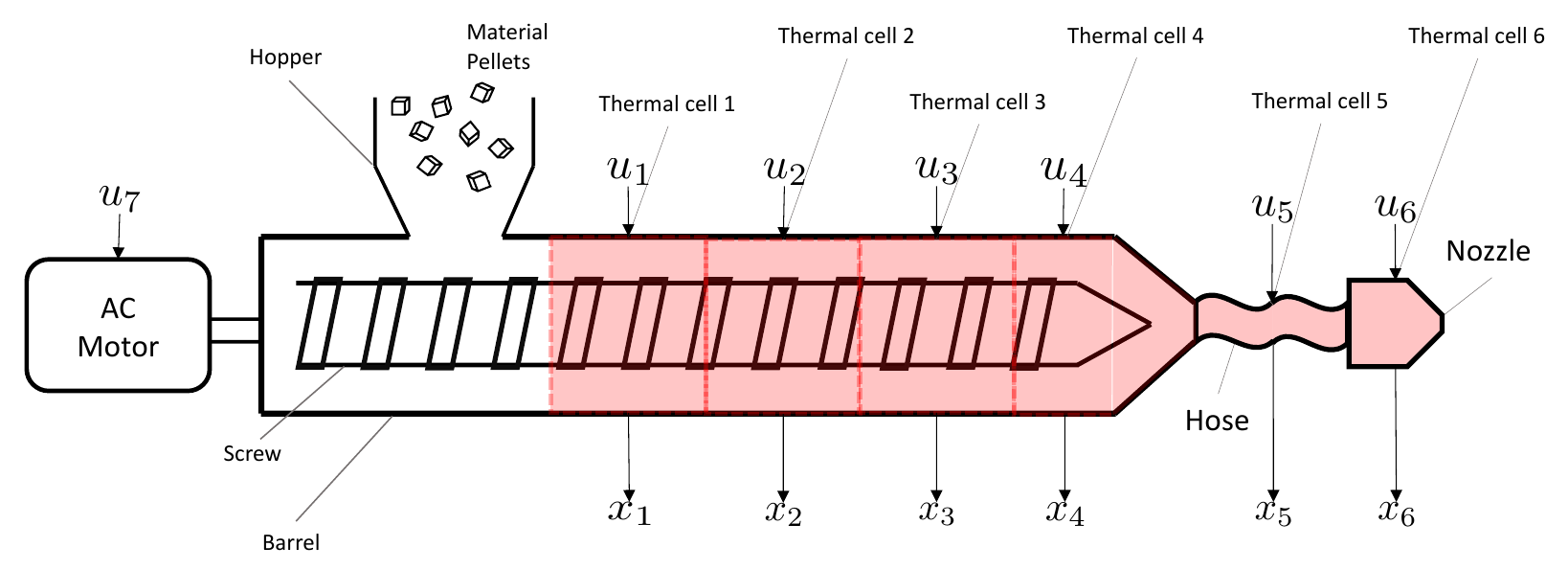}
    \caption{Extruder setup for BAAM system. Inputs 1 through 6 represent the heat applied to each thermal cell and $u_7$ represents the input from the motor. There are four heating zones inside the barrel, one on the hose and one on the nozzle. The red dotted lines denote the heating zones in the barrel. The system states represent the temperatures in each thermal cell.}
    \label{model schematic}
\end{figure}

We focus on the heating and melting of the material, specifically PLA during an MEX process. To that end, we will be using a model introduced in \cite{gootjes2017applying}. It is a discrete-time linear model in state space form and was obtained using system identification methods using input and response data from a BAAM system. Specifically, we are dealing with a model of the temperatures within the extruder of the system. The extruder comprises five key components, namely the hopper, the screw, the barrel, the hose, and the nozzle, with thermal energy being supplied by four heaters in the barrel one in the hose and one in the nozzle. An AC motor powers the rotation of the screw, which affects the temperature in each thermal cell as well as the adjacent cells. Figure \ref{model schematic} shows the heating system within the extruder, as well as indicates where the inputs are applied.  The dynamics of the system can be expressed in standard state space form as:

\begin{equation}\label{sys1}
    x(t+1)=Ax(t)+Bu(t), \; t \geq t_0
\end{equation}
where $x(t)\in \mathbb{R}^n$ is the system state vector at time $t$ and $u(t)\in \mathbb{R}^m$ is the system input at time $t$. $A\in \mathbb{R}^{n\times n}$ and $B\in \mathbb{R}^{n \times m}$ are the state and input matrices respectively. In our case the extruder temperature model, $x(t)\in \mathbb{R}^6$ is the set of the temperatures in each thermal cell and $u(t)\in \mathbb{R}^7$ is the input provided by the heaters and the motor for each time point $t$. 

The control objective is defined as the tracking of a reference signal $r(t) \in \mathbb{R}^n$ and can be expressed mathematically as in \cite{anderson2007optimal} through the performance index 
\begin{equation}\label{perf_index1}
    V_f(x,r,u)= \sum_{t=t_0+1}^T \lbrace \lbrack x(t)-r(t)\rbrack^T Q\lbrack x(t)-r(t)\rbrack + u^T(t-1)Ru(t-1)  \rbrace
\end{equation}
which quantifies how well the reference signal is tracked under input $u$, starting from time $t_0$ until the end of a time horizon $T$. $Q$ and $R$ are matrices applying appropriate weights to the tracking error and the control function respectively. In certain applications there is either no specifically defined time horizon, or it is too large. In those cases the infinite horizon performance index can be defined as in \cite{kiumarsi2014reinforcement}
\begin{equation}\label{perf_index2}
    V_{inf}(x,r,u)= \sum_{t=t_0+1}^{\infty} \gamma^{t-t_0-1}\lbrace \lbrack x(t)-r(t)\rbrack^T Q\lbrack x(t)-r(t)\rbrack + u^T(t)Ru(t)  \rbrace
\end{equation}
where $0< \gamma \leq 1$ is the discount factor, whose role is to emphasise short-term costs as opposed to costs in the more distant future. In the RL literature the performance index is often called the value function. 

The optimal tracking problem is the defined as the search for the optimal control $u^*$ such that the state $x$ tracks the reference $r$ while minimising the performance index. In the next sections, we explore solutions to the optimal tracking problem for the finite and infinite horizon cases. We consider the ideal scenario of having access to the system dynamics as well as the more realistic scenario of an unknown model. In the latter case, the problem is solved strictly using data.

\section{Model-based Control} \label{section model based}
The main idea in model-based control is the design of a controller offline, using knowledge from the model that describes the system. The controller is then brought online, interacting with the system, using the state of the system to adjust its input in order to achieve the desired behaviour. 

\subsection{Finite Horizon}
Assume the finite horizon, discrete, linear time-invariant (LTI) system defined in \eqref{sys1} and the corresponding finite horizon performance index \eqref{perf_index1}. The solution to the optimal tracking problem, as thoroughly studied in \cite{anderson2007optimal} can be found as follows:

\begin{equation}\label{lqt}
    u^*(t) = - (B^T S(t+1) B + R)^{-1} B^T \lbrack S(t+1) A x(t) + b(t+1)-Q r(t+1) \rbrack 
\end{equation}
for all $ t_0\leq t \leq T-1$, where

\begin{eqnarray}\label{b_lqt}
b(t)&=& (A^T+KB^T)(b(t+1)-Qr(t+1)), \; \text{ with } b(T)=0 \text{ and } \\
K^T &=& - ( B^T S(t+1) B +R )^{-1}B^T S(t+1)A 
\end{eqnarray}
and
\begin{eqnarray}\label{s_lqt}
S(t) &=& A^T \lbrace S(t+1)-S(t+1)B ( B^T S(t+1) B +R )^{-1}B^T S(t+1)\rbrace A + Q , \nonumber \\ & & \; S(T)=Q
\end{eqnarray}

In order to obtain a formula for $u^*(t)$ for all $ t_0\leq t \leq T-1$, given system \eqref{sys1} and performance index \eqref{perf_index1}, Equations \eqref{b_lqt} and \eqref{s_lqt} need to be solved offline, backwards in time to determine $b(t)$ and $S(t)$ for all $t$. Then the control $u^*$ can be used online, through Equation \eqref{lqt}, for all $t$. A step by step description of the solution to the finite-horizon LQT problem is given in Algorithm \ref{lqt_alg}.

\begin{algorithm}
		\SetAlgoLined
		\KwResult{Optimal controller $u^*(t)$ and optimal trajectory $x^*(t)$}
		\KwIn{State and input matrices $A,B$, Weighting matrices for performance index $Q,R$, Reference trajectory $r(t)$, initial state $x(t_0)$, Time horizon $T$}
		\begin{enumerate}
		    \item Initialise $b(T)=0$ and $S(T)=Q$.
		    \item Solve for all $b(t)$ and $S(t)$ backwards in time (offline):
		\begin{enumerate}[(i)]
		    \item Set $t=T-1$.
		    \item $K^T = - ( B^T S(t+1) B +R )^{-1}B^T S(t+1)A $
		    \item $b(t)= (A^T+KB^T)(b(t+1)-Qr(t+1))$
		    \item $S(t) = A^T \lbrack S(t+1)-S(t+1)B ( B^T S(t+1) B +R )^{-1}$ $  B^T S(t+1)\rbrack A+ Q $
		    \item Set $t=t-1$.
		\end{enumerate}
		Repeat steps $(ii)-(v)$ while $t\geq t_0+1$.
		\item Determine and apply optimal controller to system (online):
		\begin{enumerate}[(a)]
		    \item Set $t=t_0$
		    \item $u^*(t) = - (B^T S(t+1) B + R)^{-1} B^T \lbrack S(t+1) A x(t) + b(t+1)$ $-Q r(t+1) \rbrack $
		    \item $x^*(t+1)=Ax^*(t)+Bu^*(t)$
		    \item Set $t=t+1$.
		\end{enumerate}
		Repeat steps $(b)-(d)$ while $t\leq T-1$.
		\end{enumerate}
		\caption{ \label{lqt_alg} Finite-Horizon LQT}
\end{algorithm}

An important prerequisite for the design of the above controller is the predefined time horizon as it is necessary for the calculation of the terms $b(t)$ and $S(t)$. However, this approach is not feasible in cases where the horizon is not explicit. This lead to the development of infinite horizon approaches.

\subsection{Infinite Horizon}
 The infinite Horizon Optimal Tracking Problem can be addressed in a slightly different way to the previous case. Specifically, the performance index is now of the form \eqref{perf_index2}. Assume that the reference signal $r(t) \in \mathbb{R}^n$ is generated by 
\begin{equation}\label{reference_sys}
    r(t+1)=Fr(t).
\end{equation}
Consider the augmented state, including the system state and reference,
\begin{equation*}
    X(t)=\begin{bmatrix} x(t) \\ r(t)
    \end{bmatrix}
\end{equation*}
and construct the augmented system state equation
\begin{equation}
    X(t+1)= \begin{bmatrix} x(t+1) \\ r(t+1) 
    \end{bmatrix}= \begin{bmatrix} A& 0 \\ 0&F
    \end{bmatrix}\begin{bmatrix} x(t) \\ r(t)
    \end{bmatrix}+\begin{bmatrix} B \\ 0
    \end{bmatrix}u(t) = TX(t)+B_1u(t).
\end{equation}
The performance index \eqref{perf_index2} can be written in terms of the augmented state as
\begin{equation}\label{perf_index3}
    V_{inf}(X,u)= \sum_{t=k}^{\infty} \gamma^{t-k}\lbrace X^T(t) Q_1X(t) + u^T(t)Ru(t)  \rbrace
\end{equation}
where $Q_1=  \begin{bmatrix} I & -I \end{bmatrix}^TQ \begin{bmatrix} I & -I \end{bmatrix}$.
The solution to the optimal tracking problem is a policy of the form \cite{kiumarsi2014reinforcement}:
\begin{equation}\label{policy_lqt_rl}
    u(t)=-K X(t)
\end{equation}
where \begin{equation}\label{gain_lqt_rl}
    K= (R+\gamma B_1^TPB_1)^{-1}\gamma B_1^TPT
\end{equation}
and $P$ is the solution to the Algebraic Riccati Equation (ARE):
\begin{equation}\label{riccati}
    Q_1 - P +\gamma T^TPT - \gamma^2T^TPB_1 (R+\gamma B_1^TPB_1)^{-1}B_1^TPT=0.
\end{equation}

An alternative to solving the ARE can be obtained with through Reinforcement Learning (RL). The Lemma proven in \cite{kiumarsi2014reinforcement,Lewis}, states that the value function \eqref{perf_index3} can be written in quadratic form 
\begin{equation}\label{quad_form_value_function}
    V(x,r,u)= V(X)=\frac{1}{2}X^T(t)PX(t)
\end{equation}
for some matrix $P=P^T>0$, any stabilising policy \eqref{policy_lqt_rl} and reference signal \eqref{reference_sys}. This form of the value function gives rise to the LQT Bellman equation 
\begin{eqnarray}\label{bellman}
    X^T(t)PX(t)&=&X^T(t)Q_1X(t)+u^T(t)Ru(t) \nonumber\\
    & &+\gamma X^T(t+1)PX(t+1).
\end{eqnarray}
Consider a specific stabilizing policy $K^*$ in \eqref{bellman}. That gives rise to the Lyapunov equation
\begin{equation}\label{lyapunov}
    P=Q_1 +(K^*)^TRK^*+\gamma (T+B_1K^*)^TP(T+B_1K^*).
\end{equation}
Iteratively solving the Lyapunov equation while updating the optimal control estimate can effectively lead to the solution to the Optimal Tracking Problem.  

A step by step description of the solution to the infinite-horizon problem can be found in Algorithm \ref{model_based_ql_alg}.
\begin{algorithm}
		\SetAlgoLined
		\KwResult{Optimal control gain $K^*$, optimal controller $u^*(t)$ and optimal trajectory $x^*(t)$}
		\KwIn{State and input matrices $A,B$, Weighting matrices for performance index $Q,R$, Initial state $x(t_0)$, Initial reference signal $r(t_0)$,  Discount factor $\gamma$, Initial control policy $K^0$, Error threshold $\epsilon$}
		\begin{enumerate}
		    \item Create augmented system matrices $T=\begin{bmatrix}
	A & 0 \\ 0 & F
	\end{bmatrix},B_1=\begin{bmatrix}
	B \\ 0 
	\end{bmatrix}$.
	\item Augment weighting matrix $Q_1=  \begin{bmatrix} I & -I \end{bmatrix}^TQ \begin{bmatrix} I & -I \end{bmatrix}$.
	\item Determine optimal control gain $K^*$:
	\begin{enumerate}[(i)]
	    \item Set $j=0$.
	    \item Solve the LQT Lyapunov equation $P^{j+1}=Q_1+(K^j)^T RK^j +\gamma (T-B_1K^j)^TP^{j+1}(T-B_1K^j)$.
	\item Compute the control gain estimate 
		$K^{j+1}= (R+\gamma B_1^T P^{j+1} B_1)^{-1}$  $\gamma B_1^TP^{j+1}T$.
		\item Set $j=j+1$.
	\end{enumerate}
	Repeat steps $(ii)-(iv)$ while $\vert K^{j}-K^{j-1}\vert>\epsilon$. 
	\item Set $K^*=K^j$.
	\item Determine and apply optimal controller to system (online):
		\begin{enumerate}[(a)]
		    \item Design initial augmented state $X(t_0)=\begin{bmatrix} x(t_0) \\ r(t_0)\end{bmatrix}$
		    \item Set $t=t_0$.
		    \item $u^*(t) = - K^* X(t) $
		    \item $X^*(t+1)=TX^*(t)+B_1u^*(t)$
		    \item Retrieve optimal state $x^*(t)$.
		    \item Set $t=t+1$.
		\end{enumerate}
		Repeat steps $(c)-(e)$ for the desired length of the experiment or simulation.
		\end{enumerate}
		\caption{ \label{model_based_ql_alg} Infinite-Horizon LQT with RL}
	\end{algorithm}   

\section{Data-driven Q-Learning Control}  \label{section q-learning}
In the previous sections, the Linear Quadratic Tracking problem was thoroughly studied in the case where a model of the system of interest is available. However, this is not the case in most applications. In most instances, a model is either entirely unknown or not accurate enough to be considered ground truth. This is the motivation behind the constantly developing field of data-driven control. In this section we explore Reinforcement Learning based algorithms, specifically ones in the Q-learning \cite{watkins1989learning} family, which are named after the use of Q-functions as value functions. 

We explore an online Q-learning solution to the LQT problem that is dependent solely on acquired data of the state trajectories \cite{kiumarsi2014reinforcement}. The data is used to estimate the solution to the algebraic Riccati equation \eqref{riccati} following the analysis described below. 

Let our system be described by the augmented state equation of the previous section \begin{equation}\label{augmented_state}
    X(t+1)=TX(t)+B_1u(t)
\end{equation}
and let the weighting matrices $Q_1$ and $R$ be defined as in \eqref{perf_index3}
The Q-function is defined with the help of the Bellman equation as the sum of the current cost and the discounted future cost. 
\begin{equation}\label{q_func_bellman}
    Q(x,r,u)=\frac{1}{2}X^T(t)Q_1X(t) + \frac{1}{2}u^T(t)Ru(t) + \frac{1}{2}\gamma X^T(t+1)PX(t+1)
\end{equation}
where P is the solution to \eqref{riccati}. Using the system dynamics, the Q-function can be written as:
\begin{eqnarray}
    Q(X,u)&=&\frac{1}{2}X^T(t)Q_1X(t) + \frac{1}{2}u^T(t)Ru(t) \\ 
    & &+ \frac{1}{2}\gamma (TX(t)+B_1u(t))^TP(TX(t)+B_1u(t)) \\
    &=& \frac{1}{2}\begin{bmatrix} X(t) \\ u(t)
    \end{bmatrix}^T \begin{bmatrix} Q_1+\gamma T^TPT & \gamma T^TPB_1 \\ \gamma B_1^TPT &R+\gamma B_1^TPB_1
    \end{bmatrix}\begin{bmatrix} X(t) \\ u(t)
    \end{bmatrix}\\
    &=& \frac{1}{2}\begin{bmatrix} X(t) \\ u(t)
    \end{bmatrix}^T H \begin{bmatrix} X(t) \\ u(t)
    \end{bmatrix}
\end{eqnarray}
where H is the symmetric kernel matrix. H can be written in block form as:

\begin{equation}
    H=\begin{bmatrix}
    H_{XX} & H_{Xu} \\ H_{uX} & H_{uu}
    \end{bmatrix}
\end{equation}

Due to the quadratic form of the Q-function, its minimisation can be reduced to the solution of $\frac{\partial Q(X(t),u(t))}{\partial u(t)}=0$. This yields
\begin{equation}
    u(t)=-H_{uu}^{-1}H_{uX} X(t) = -(R+\gamma B_1^TPB_1)^{-1}\gamma B_1^TPTX(t)
\end{equation}
which is the same result as in equations \eqref{policy_lqt_rl} and \eqref{gain_lqt_rl}. 
It is evident that the key to obtaining the optimal control policy in this formulation is the kernel matrix $H$. This matrix is what will be estimated through the Q-learning algorithm.

Setting $Z(t)=\begin{bmatrix} X(t) \\ u(t) \end{bmatrix}$, the Q-function becomes
\begin{equation}
    Q(Z)= \frac{1}{2}Z^T(t)HZ(t)
\end{equation}
and consequently, \eqref{q_func_bellman} can be written as
\begin{equation}\label{final_bellman}
    Z^T(t)HZ(t)=X^T(t)Q_1X(t) + u^T(t)Ru(t) +Z^T(t+1)HZ(t+1).
\end{equation}
In the above expression, the only element that requires knowledge of the system dynamics is $H$. However, it can be approximated using obtained data, namely a set of state trajectories with the corresponding references and inputs $(x,r,u)$. The weighting matrices $Q_1,R$ are predefined by the design of the optimisation problem. In order to estimate $H$, it needs to be isolated, meaning that the quadratic form of \eqref{final_bellman} needs to be transformed into an equivalent form. We start by vectorising \eqref{final_bellman} and then use the "vector trick" associated with the Kronecker product.
\begin{eqnarray}
    vec \left( Z^T(t)HZ(t) \right) &=&  X^T(t)Q_1X(t)  +  u^T(t)Ru(t) \nonumber\\ & &+ \gamma vec \left( Z^T(t+1)HZ(t+1) \right)\\
    \left( Z^T(t)\otimes Z^T(t) \right)vec(H)&=& X^T(t)Q_1X(t)  +  u^T(t)Ru(t)\nonumber\\ & &+ \gamma  \left(Z^T(t+1)\otimes Z^T(t+1) \right) vec(H) \label{final_value_func}
\end{eqnarray}

Equation \eqref{final_value_func} can be iteratively solved through the \textit{Value Iteration} algorithm \cite{sutton2018reinforcement} which consists of two steps:

\begin{enumerate}
    \item \textbf{Policy Evaluation} 
    \begin{eqnarray}
        \left( Z^T(t)\otimes Z^T(t) \right)vec(H^{i+1})&=& X^T(t)Q_1X(t)  +  u^T(t)R u(t)\label{policy_eval_state} \\ & &+ \gamma  \left(\hat{Z}^T(t+1)\otimes \hat{Z}^T(t+1) \right) vec(H^i)\nonumber
    \end{eqnarray}
    where $\hat{Z}(t+1)=\left\lbrace X(t+1),u^i(t+1)\right\rbrace$
    \item \textbf{Policy Improvement}
    \begin{eqnarray}
        u^{i+1}(t) &=& -(H^{i+1}_{uu})^{-1}H^{i+1}_{uX} X(t)
    \end{eqnarray}
\end{enumerate}
The Policy evaluation step can be solved through the Least Squares algorithm using measured data $Z(t),Z(t+1)$ and calculating the cost term $X^T(t)Q_1X(t)  +  (u^i(t))^TR u^i(t)$ for each data point. Matrix $H$ is an $(2n+m)\times (2n+m+1)$ symmetric matrix which means that its determination is a problem with $(2n+m)\times (2n+m+1)/2$ degrees of freedom. Hence, at least $(2n+m)\times (2n+m+1)/2$ data points are required to solve \eqref{policy_eval_state}. It is worth noting that in practice, many more data points are usually needed, especially when dealing with larger state and action spaces. Another important note is the need for including a regularisation step in the Least Squares algorithm. When the data are highly correlated, the inversion step becomes numerically unstable. This problem can be reduced with the choice of an appropriate regularisation parameter $\mu$. A step-by-step description of the solution to the LQT problem using Q-learning for the state feedback case can be found in Algorithm \eqref{alg_ql_state}.

\begin{algorithm}
		\SetAlgoLined
		\KwResult{Kernel matrix $H$, optimal controller $u^*(t)$ and optimal trajectory $x^*(t)$}
		\KwIn{Weighting matrices for performance index $Q_1,R$, regularisation parameter $\mu$, Discount factor $\gamma$, initial kernel matrix estimate $H^0$}
		\begin{enumerate}
		    \item Collect $N>>(2n+m)\times (2n+m+1)/2$ consecutive data tuples \\$Z(t)=\left\lbrace X(t),u(t)\right\rbrace$ where $X(t)=\left\lbrace x(t),r(t)\right\rbrace $.
		    \item For each $Z(t)$, calculate the Kronecker product $K_Z(t):=Z^T(t)\otimes Z^T(t)$.
		    \item Set $i=0$.
		        \item For t=1 to N:
		        \begin{enumerate}[(i)]
		            \item Determine $u^i(t+1)= -(H^{i}_{uu})^{-1}H^{i}_{uX} X(t+1)$
		        \item Create the updated augmented state $\hat{Z}(t+1)=\left\lbrace X(t+1),u^i(t+1)\right\rbrace$ and the corresponding Kronecker product $K_{\hat{Z}}(t+1)=\hat{Z}^T(t+1)\otimes \hat{Z}^T(t+1)$
		        \item Compute the future expected cost term $c(t)=X^T(t)Q_1X(t) $ $ +  u^T(t)R u(t)+ \gamma  K_{\hat{Z}}(t+1) vec(H^i)$
		        \end{enumerate}
		        \item Calculate the sums $L=\sum_1^N K^T_Z(t)K_Z(t)$ and $R=\sum_1^N K^T_Z(t)c(t)$ .
		        \item Regularise matrix L: $L_r=L+\mu I$
		        \item Produce new H estimate $H^{i+1}=L_r^{-1}R$
		    \item i=i+1
		\end{enumerate}
		Repeat steps (4)-(8) until convergence of the estimate of $H$ or for a predetermined number of iterations.
		\caption{ \label{alg_ql_state} State Feedback Q-Learning}
\end{algorithm}

\section{Simulations and Numerical Results} \label{section results}
The model used throughout the simulations is identified in \cite{gootjes2017applying}. The system matrices $A\in \mathbb{R}^{6\times 6}$ and $B\in \mathbb{R}^{6\times 7}$ were determined through system identification approaches and were found to be:

\begin{eqnarray}
    A&=& \begin{bmatrix}
    0.992 & 0.0018 & 0 & 0 & 0 & 0 \\ 
    0.0023 & 0.9919 & 0.0043 & 0 & 0 & 0 \\ 
    0 & -0.0042 & 1.0009 & 0.0024 & 0 & 0 \\ 
    0 & 0 & 0.0013 & 0.9979 &0&0 \\ 
    0&0&0&0&0.9972&0 \\
    0&0&0&0&0&0.9953
    \end{bmatrix} \\ 
    B&=& \begin{bmatrix}
    1.0033 & 0 & 0 & 0 & 0 & 0&-0.2175 \\ 
    0 & 1.0460 & 0 & 0 & 0 & 0&-0.0788 \\ 
    0 & 0 & 1.0326 & 0 & 0 & 0 &-0.0020\\ 
    0 & 0 & 0 & 0.4798 &0&0 &-0.0669\\ 
    0&0&0&0&0.8882&0 &0.1273\\
    0&0&0&0&0&1.1699&-0.1792
    \end{bmatrix} 
\end{eqnarray}

This model introduced will be used for two main purposes. Firstly, it will be the basis around which the model-based control algorithms will be built. Secondly, it will be used throughout the simulation steps, generating data with which to train the data-driven algorithms. 

We set up an optimisation problem to be solved with all methods mentioned in the above sections. The goal of the optimisation is to control the behaviour of the system introduced in Section \ref{section_model} so as to maintain the temperatures of all system states at some predefined set of increasing values, similar to Figure \ref{temp example}. We chose the reference trajectory to be the vector of increasing values: $r(t)=\left\lbrace 155,160,165,170,180,190\right\rbrace$\textdegree C $\in \mathbb{R}^6$ for all time points $t$, so that the last two heaters bring PLA within its melting temperature range. This also implies that the reference generator matrix $F=I_6$.

\begin{figure}[ht]
	\begin{subfigure}{.5\textwidth}
		\includegraphics[width=\linewidth,height=4cm]{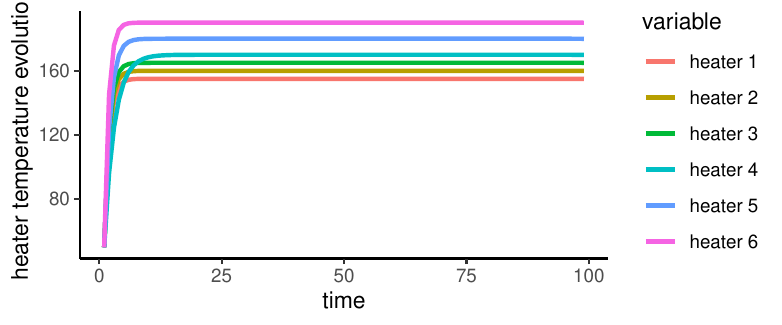}
		\caption{\label{fig:lqt anderson}Finite Horizon LQT}
	\end{subfigure}%
	\begin{subfigure}{.5\textwidth}
		\includegraphics[width=\linewidth,height=4cm]{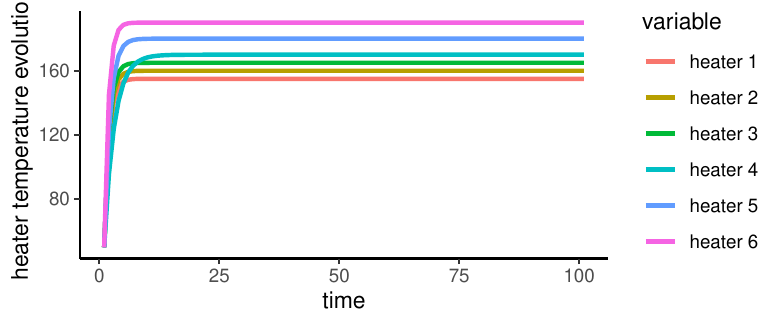}
		\caption{\label{fig:lqt infinite}Infinite Horizon LQT}
	\end{subfigure}
	\caption{\label{Model based controllers}Model based control. Both algorithms successfully bring and maintain the heater temperatures within a small margin of the optimal values. In the finite horizon case, the convergence more quickly, while in the infinite horizon case, more time steps are required before the temperatures converge to the optimal.}
	\end{figure}

Starting with the model-based methods, we initially focus on the finite-horizon LQT problem. The parameters that need to be determined for the implementation of Algorithm \ref{lqt_alg} are the weighting matrices $Q$ and $R$, the initial state $x(t_0)$ and the time horizon $T$. For $Q,R$ we chose identity matrices of appropriate dimensions, namely $Q=I_6$ and $R=I_7$. The initial state was arbitrarily chosen to be $x(t_0) =\left\lbrace 50,\dots , 50\right\rbrace \in \mathbb{R}^6$ and the time horizon was chosen to be $T=100$. Being the most reliable and efficient controller, implementation of Algorithm \ref{lqt_alg} results in an ideal system behaviour, where the optimal is achieved and maintained within very few steps. Specifically, all 6 heaters achieve temperatures within less than $0.09 $\textdegree C of the optimal temperature vector, within 16 steps. The values that the heater temperatures converge to are $x^*=\left\lbrace154.918,159.985, 164.992, 169.922, 179.926, 189.955\right\rbrace$. This result can be seen in Figure \ref{fig:lqt anderson}. 

Next, we explore the case of the infinite time horizon, while still assuming full knowledge of the system dynamics, making use of Algorithm \ref{model_based_ql_alg}. The weighting matrices for the performance index $Q,R$ and the initial state $x(t_0)$ are chosen to be the same as the finite horizon case. There are a few additional parameters that need to be determined. The discount factor is chosen as $\gamma=0.99$, and the initial control policy $K^0$ is chosen as a $7\times12$ matrix populated by normal random numbers between with mean 0 and standard deviation 10. Finally, the error threshold needed for the estimation of the control policy $K$ is chosen to be $\epsilon=0.1$. The RL-based controller is successful in bringing and maintaining all temperatures within 0.04 degrees of the optimal and more specifically, the temperature vector converges to the vector $x^*=\left\lbrace154.972, 159.996, 164.966, 169.989, 179.997, 189.955\right\rbrace$. 17 time steps are needed for the heaters to approach within $0.1$ \textdegree C of the optimal and finally settle to their final values after 36 steps. This result is visualised in Figure \ref{fig:lqt infinite} where we plot the state trajectories for the first 100 time steps. 

Table \ref{table_model_based} summarises all the parameters used in the model-based approaches to produce Figure \ref{Model based controllers}. 
	
	\begin{table}[ht]
	\centering
	\begin{tabular}{|c|c||c|} 
	\hline
	\multicolumn{3}{|c|}{Table of parameters - Model based methods}\\
		\hline
		\hline
		Parameter&Symbol&Value\\
		\hline
		\multicolumn{3}{|c|}{Universal parameters}\\
		\hline 
		Tracking reference & r & $\left\lbrace 190,180,170,165,160,155 \right\rbrace$ \\
		Initial state & $x(t_0)$ & $\left\lbrace 50, \dots, 50 \right\rbrace$\\
		Tracking error weighting matrix&Q & $I_6$\\ 
		Input weighting matrix&R & $I_7$\\ 
		\hline
		\multicolumn{3}{|c|}{Finite Horizon LQT}\\
		\hline 
		Time horizon & T & 100  \\ 
		\hline
		\multicolumn{3}{|c|}{Infinite Horizon LQT}\\
		\hline 
		Discount factor & $\gamma$& 0.99\\ 
		Error threshold & $\epsilon$& 0.1 \\ 
		\hline
	\end{tabular} 
	\caption{\label{table_model_based}Parameters chosen for model-based control implementations. The information listed includes the parameter names, the corresponding symbol and the value that was chosen grouped by the algorithms in which they were used.}
\end{table}

Finally, we focus on the Q-learning algorithm introduced in Section \ref{section q-learning}. As this is a learning algorithm, there is no access to the model of the system for the design of the controller, i.e., matrices $A,B$ are unknown. Instead, the model introduced in Section \ref{section_model} is used only as a data generator. In a real application of this learning approach, the data will be made available from sensorising the process itself, meaning that the steps regarding the data generation will not be necessary. 
	
The simulated data needs to be varied enough so that it holds sufficient information about the system. This can be achieved by applying an input that is persistently exciting \cite{yuan1977probing,de2019formulas}, which can be achieved by the inclusion of probing noise. For our application, we consider the following random vectors 
	
 \begin{equation}
    \omega_1 \sim \mathcal{N}(0,\sqrt{\sigma}I_7), \omega_2, \dots \omega_8 \sim \mathcal{N}(0, I_7)
\end{equation}
where $\sigma=0.5$. We design the probing noise vector for time $t$ by summing sinusoidal signals of various frequencies and periods, as follows:
	
\begin{eqnarray}\label{noise_vector}
    \omega_{pr}(t) &= &\omega_1 + 10 \sin \left( \frac{\omega_2\pi t}{5}\right)+ 8 \sin \left( \frac{2\omega_3\pi t}{5}\right) +  7 \sin \left( \frac{3\omega_4\pi t}{5}\right) \nonumber \\ 
    & &+ 6 \sin \left( \frac{4\omega_5\pi t}{5}\right) +  4 \sin \left( \frac{5\omega_6\pi t}{5}\right) +  3 \sin \left( \frac{6\omega_7\pi t}{5}\right)  \nonumber \\ 
    & &+  0.5 \sin \left( \frac{7\omega_8\pi t}{5}\right)
\end{eqnarray}
 
It is important that the process remains stable during the simulations so that the states do not end up diverging towards infinity. This means that the input, apart from being persistently exciting, also needs to be stabilising. To that end, we use the following stabilising gain: 
\begin{equation}
	K=
	    \begin{bmatrix*}[r]
	    0.7395 & -0.0076 & -0.0003 & -0.0264 & 0.0194 & -0.0170 \\
	    -0.0076 & 0.7430 & 0.0031 &-0.0093 & 0.0068 & -0.0060 \\ 
	    -0.0003 & -0.0033 & 0.7599 & 0.0021 & 0.0002 & -0.0002 \\
	    -0.0126 & -0.0042 & 0.0016 & 1.0971 & 0.0092 & -0.0079 \\
	    0.0171 & 0.0058 & 0.0002 & 0.0170 & 0.8179 & 0.0108 \\
	    -0.0198 & -0.0067 & -0.0002 & -0.0193 & 0.0143 & 0.6823 \\
	    -0.1525 & -0.0519 & -0.0018 & -0.1412 & 0.1091 & -0.0977
	    \end{bmatrix*}
\end{equation}

Considering all the above, the input to the system for time $t$ can be defined as
\begin{equation}
    u(t) = -Kx(t) + \omega_{pr}(t)
\end{equation}
where $\omega_{pr}$ is the probing noise vector, as defined in \eqref{noise_vector}. Starting with the same state initialisation as the previous cases, namely $x(t_0)=\left\lbrace50,\dots,50 \right\rbrace$ and initialising the input function as $u(t_0)=\omega_1$, we produce a set of 2000 data points that will be used for training. 

In order to train the Q-learning based controller, a few more parameters need to be determined. Firstly, the weighting parameters of the value function are chosen as above, $Q=I_6$ and $R=I_7$. The discount factor was chosen to be $\gamma=0.99$ and the regularisation parameter $\mu=0.001$. Finally, the initial estimate of the kernel matrix is chosen to be $H^0=I_{19}$. We train the algorithm until the error between the estimates of $H$ becomes smaller than $0.001$ or for a maximum of 30 iterations. Table \ref{table_data_driven} summarises all parameter values used for the implementation of Q-learning for the state feedback problem. 

\begin{figure}[ht]
\centering
		\includegraphics[width=0.5\linewidth,height=4cm]{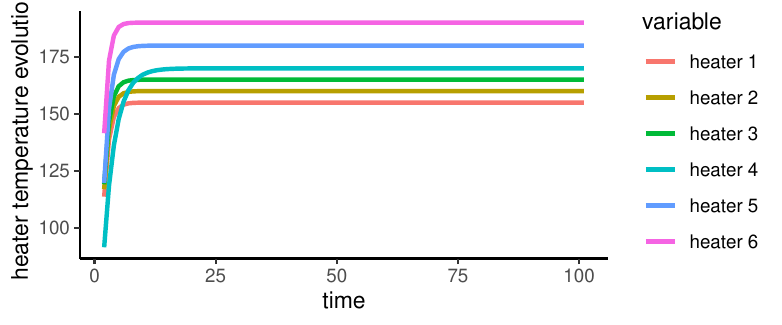}
		\caption{\label{fig:q-learning}State Feedback Q-learning. The learning algorithm obtains a control gain that is sub-optimal but successfully brings and maintains all temperatures within 0.11 degrees of the reference.}
	\end{figure}

\begin{table}[ht]
	\centering
	\begin{tabular}{|c|c||c|} 
	\hline
	\multicolumn{3}{|c|}{Table of parameters - Q-Learning}\\
		\hline
		\hline
		Parameter&Symbol&Value\\
		\hline
		Tracking reference & r & $\left\lbrace 190,180,170,165,160,155 \right\rbrace$ \\
		Initial state & $x(t_0)$ & $\left\lbrace 50, \dots, 50 \right\rbrace$\\
		Tracking error weighting matrix&Q & $I_6$\\ 
		Input weighting matrix&R & $I_7$\\ 
		Discount factor & $\gamma$ & 0.99 \\
		Regularisation parameter & $\mu$ & 0.0001 \\ 
		Initial kernel matrix estimate &$H^0$& $I_{19}$\\
		\hline
	\end{tabular} 
	\caption{\label{table_data_driven}Parameters chosen for the implementation of Q-learning. The information listed includes the parameter names, the corresponding symbol and the value that was chosen.}
\end{table}
\sloppy
Figure \ref{fig:q-learning} displays the state trajectories that result by applying to the system the controller learned with the use of Algorithm \ref{alg_ql_state}. The trajectories are plotted over 100 time steps. The states converge to the vector $x^*=\left\lbrace 154.895, 159.978, 164.990, 169.939, 179.922 189.929 \right\rbrace$ which means they achieve temperatures within $0.11 $\textdegree C of the optimal. The state trajectories arrive within that error margin of the optimal in 20 time steps and convergence to their final values is completed within 43 steps. The results have a slight deviation from the optimal results produced by the model-based methods, but are close enough that plots \ref{fig:lqt infinite} and \ref{fig:q-learning} seem almost identical to the naked eye. This implies that the Q-learning algorithm converges to some sub-optimal kernel matrix H (and corresponding gain K). In fact, the Frobenious norm of the absolute difference between the optimal gain $K^*$ as determined by Algorithm \ref{model_based_ql_alg} and the estimated gain $\hat{K}=H_{uu}^{-1}H_{uX}$ from Algorithm \ref{alg_ql_state} is 0.756. The mean absolute error of the individual elements populating the gain matrix is  0.0803.

	A different choice in the discount factor $\gamma$, would result in a different tracking performance. For example, choosing $\gamma=0.95$ would cause the state trajectories to arrive and remain within $0.07$\textdegree C of the optimal, and choosing $\gamma=0.9$ would bring the temperatures within $0.14$\textdegree C of the reference trajectory.

\begin{figure}[ht]
\centering
		\includegraphics[width=0.5\linewidth,height=4cm]{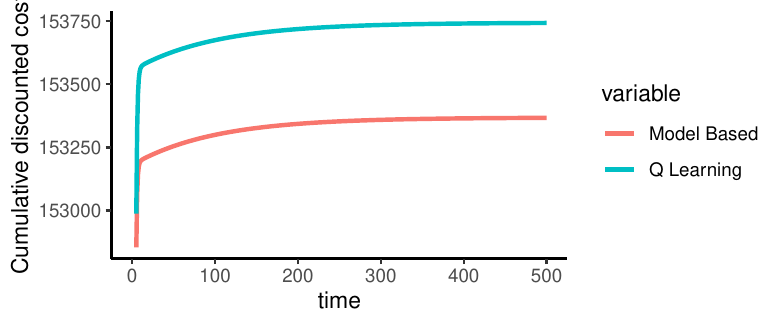}
		\caption{\label{fig:cost comparison}Comparison of the cumulative costs between the model-based and model-free algorithms. When the system dynamics are known, the performance index is more effectively minimised.}
	\end{figure}

Another way to compare the performance of the Q-learning data-driven controller with its model-based counterpart is to look into the minimisation of the performance index \eqref{perf_index3}. This is a particularly important measure of performance, considering that the data-driven controller seemingly performs better at tracking the optimal vector. However, for both algorithms, the optimisation goal is the minimisation of the infinite sum of discounted costs, which includes both the tracking cost and the cost associated with the input. Specifically, the model-based algorithm after 500 steps achieves a cumulative cost of $153,367 $ while the data-driven method reaches the value $153,743$. In around 900 steps both algorithms have achieved their final performance index which is $153,368 $ and $  153,745$ for the model-based and model-free cases respectively. Figure \ref{fig:cost comparison} shows the two cumulative costs over the first 500 time steps.
The better tracking performance obtained by the model-based algorithm is verified by looking at the part of the cost corresponding to the reference tracking. For the model-based case, the total tracking cost is $107,830$ while for the data-driven algorithm this value is $111,676$. However, when looking at the cost term associated with the input, it is higher in the model based case, namely $45,538 $ compared to $42,068$ achieved by the Q-learning algorithm. This is expected, as the total cumulative costs are very very close in value. 

It is worth noting that some of the numerical results might slightly vary between different implementations, as they are dependent on the randomness introduced in the algorithm. 
  
\section{Conclusion} \label{section_conclusion}
We have studied the problem of controlling a state space system. The system at hand was the heating system within the extruder of a BAAM system. The control objective was set to be bringing and maintaining the temperatures to some predetermined set of values. The problem was approached through different lenses, namely classic Control Theory techniques as well as Reinforcement Learning approaches. Initially, the focus was on model-based methods, exploring both the finite and infinite horizon cases, the first approached through the standard LQT algorithm, the second through Reinforcement Learning. Finally we explored the case of data-driven, model-free state feedback.   

We showed that, in the case when an exact model of the system is available, the design of controllers is a deterministic process as it stems directly from the system dynamics. Both the finite and infinite horizon LQT controllers successfully achieve the optimisation goal, with a small error margin and within very few steps. 

Given that obtaining exact models of physical processes is usually not possible, with that being evident also in Additive Manufacturing, the field of Data-driven Control is being very actively researched and developed. We approached the problem of optimising the temperatures within the extruder without knowledge of its exact dynamics, using Reinforcement Learning techniques and specifically Q-learning. The assumption is made that the entire state of the system is available for measurement. With an appropriate choice of parameters and a large training data set, the obtained controller was able to achieve the desired behaviour in the system within a few time steps. This is a very encouraging result, since it implies that many systems in general and 3D printers in particular, can learn to optimise their behaviour simply by gathering data, without the need for complicated models. The more advanced data-driven controllers become, the closer the idea of truly Intelligent Manufacturing gets. 

The results in this paper are evidence to the fact that Reinforcement Learning, and more specifically Q-learning, can produce comparable results to model-based methods when controlling a state-space system. However, a few drawbacks in the approaches are highlighted. Firstly, the assumption that the entire state is exactly measurable is unrealistic. In real-life applications, measurements are obtained through sensors, which can produce errors or measure multiple states at once. Hence, output feedback approaches need to be investigated. Secondly, the training of the Q-learning algorithm requires a large amount of data which can be an issue, especially when dealing with more complicated systems. Machine learning approaches that are good at handling big data need to be explored, such as Deep Learning. These issues remain as future work.

\section*{Declarations}

\begin{itemize}
\item Funding:
This publication has emanated from research supported by a research grant from Science Foundation Ireland (SFI) under Grant Number 16/RC/3872.
\item Competing interests: The authors have no competing interests to declare that are relevant to the content of this article.
\item Ethics approval: Not applicable
\item Consent to participate: Not applicable
\item Consent for publication: Not applicable
\item Availability of data and materials: Not applicable
\item Code availability: The code used for our simulations was written in R and can be found on: 
\url{https://github.com/elenizavrakli/State-feedback-control-for-AM-temperature-model}
\item Authors' contributions:
S.D and A.P and E.Z conceived of the idea. E.Z developed the theory, performed the simulations, and wrote the manuscript with support from S.D and A.P. A.D provided the manufacturing context.
\end{itemize}

\bibliography{bibl}

\end{document}